\documentclass[oneside,11pt]{amsart}
\usepackage{amscd,amssymb}
\usepackage[all]{xy}

\title[Some remarks on the derived categories of homogeneous spaces]
 {Some remarks on the derived categories of coherent sheaves
  on homogeneous spaces}

\author{Alexander Samokhin}

\address{Institute for Information Transmission Problems, Moscow, Russia}

\email{samohin@mccme.ru}

\thanks{This work was supported in part by the French Government
  Fellowship and the RFFI award No. 02--01--22005.}

\long\def\comment#1{}

\newcommand{\Oo}{\mathcal O}
\newcommand{\Uu}{\mathcal U}
\newcommand{\Ww}{\mathcal W}
\newcommand{\Ss}{\sf S}
\newcommand{\Mm}{\sf M}
\newcommand{\Pp}{\mathbb P}

\newcommand{\Ff}{\mathcal F}
\newcommand{\Ee}{\mathcal E}
\newcommand{\Ll}{\mathcal L}

\newcommand*{\RHom}{\mathop{\mathrm RHom}\nolimits}
\newcommand*{\Hom}{\mathop{\mathrm Hom}\nolimits}
\newcommand*{\Dd}{\mathop{\mathrm D\kern0pt}\nolimits}

\newtheorem{theorem}{Theorem}[section]
\newtheorem{corollary}{Corollary}[section]
\newtheorem{lemma}{Lemma}[section]
\newtheorem{proposition}{Proposition}[section]

\newtheorem{remark}{Remark}[section]

\theoremstyle{definition}
\newtheorem{definition}{Definition}[section]

\numberwithin{equation}{section}

\long\def\comment#1{}

\comment{
\newtheorem{theorem}{Theorem}         
\newtheorem{corollary}{Corollary}
\newtheorem{lemma}{Lemma}

\newnumbered{assertion}{Assertion}    
\newnumbered{conjecture}{Conjecture}  
\newnumbered{definition}{Definition}
\newnumbered{hypothesis}{Hypothesis}
\newnumbered{note}{Note}
\newnumbered{observation}{Observation}
\newnumbered{problem}{Problem}
\newnumbered{question}{Question}
\newnumbered{algorithm}{Algorithm}
\newnumbered{example}{Example}
\newunnumbered{notation}{Notation} 

\newtheorem{remark}{Remark}[section]
}

\begin{document}

\begin{abstract}
In this paper we prove first a general theorem on semiorthogonal decompositions in
derived categories of coherent sheaves for flat families over a smooth base. We then show that the
derived categories of coherent sheaves on flag varieties of classical
type are generated by complete exceptional collections. Finally, we
find complete exceptional collections in the derived categories of
some homogeneous spaces of the symplectic groups of small rank.
\end{abstract}

\maketitle

\section{Introduction}

The study of derived categories of coherent sheaves on algebraic
varieties dates back to the late 70's when A. Beilinson described the
derived categories of projective spaces in \cite{Be}. Since then this field has experienced
a great development that put it in the forefront of modern algebraic
geometry. At the beginning of the 80's it was discovered that derived categories of
coherent sheaves appeared in a variety of areas of mathematics: e.g., in representation theory of
associative algebras, and in the study of moduli spaces, to name a
few. Other links were found later, of which the homological mirror
conjecture (\cite{Ko}), and reinterpretation of the minimal model
program in birational geometry in terms of
semiorthogonal decompositions in derived categories (\cite{BO}) are, perhaps,
the most important. Another profound connection between representations of
semisimple Lie algebras in positive characteristic and the derived categories of coherent
sheaves on Springer fibers of the Grothendieck--Springer resolutions of
nilpotent cones --- in particular, on homogeneous spaces --- was established
recently in \cite{BMR}.\par

In fact, homogeneous spaces of semisimple 
Lie groups were the first varieties studied in the context of derived categories of coherent
sheaves. Next to the paper \cite{Be} was a series of papers by
M. Kapranov who developed substantially the approach from {\it loc.cit.}.
His work culminated in the article \cite{Kap}, where the derived
categories of coherent sheaves on homogeneous spaces of Lie groups of
type ${\bf A}_n$ and on quadrics were described. Far back in the late
80's, these facts, together with other evidences, led experts to believe that the derived
category of any homogeneous space has a description similar to that
from \cite{Be} and \cite{Kap}, i.e. that such a category is generated by an
exceptional collection. More recently, the results of \cite{BLL} indicate that derived categories of coherent sheaves on 
varieties admitting pavings by affine spaces (of which
homogeneous spaces are a particular class) should have complete exceptional
collections. Such collections are a very particular
case of semiorthogonal decompositions. Over more than fifteen years, which have
passed since the appearance of \cite{Kap}, a great many powerful
techniques were invented to study semiorthogonal decompositions in derived categories.
Despite this, explicit complete exceptional collections on arbitrary homogeneous
spaces of semisimple Lie groups remain an unsolved problem
up-to-date.\par

We limit ourselves to indicating one source of our interest in these collections. 
As was shown by S. Mukai (\cite{Muk}) in his classification of prime Fano threefolds
with $b_2 = 1$, such threefolds are closely connected to
homogeneous spaces. In particular, some of these threefolds can be
obtained as linear sections of appropriate homogeneous spaces. 
It is important to construct semiorthogonal
decompositions in the derived categories of prime Fano threefolds for a number of reasons;
however, these categories are fairly complicated objects. A recent phenomenon, discovered and studied by A. Kuznetsov 
in a series of papers (\cite{Kuznetsov1}, \cite{Kuznetsov2}, \cite{Kuznetsov3}), is that
semiorthogonal decompositions in the derived categories of linear
sections of homogeneous spaces turned out
to be intrinsically related to 
those in the derived categories of linear sections of projectively
dual varieties. The approach developed in {\it loc.cit.}
permits, in particular, to describe in a uniform way the derived
categories of linear sections of some homogeneous spaces. To make it
work one needs to know the structure of 
derived category of a homogeneous space one starts with. \par

In the present paper we take some further steps towards the problem of finding complete exceptional collections
on homogeneous spaces. In Section \ref{sec:genOrlovth} we prove Theorem \ref{th:Samokh1} on
semiorthogonal decompositions for flat families over a smooth
scheme. The assumptions of Theorem \ref{th:Samokh1}, however, are
restrictive enough and, as a consequence, this theorem does not cover
such important case as flat families of quadrics. Happily, 
semiorthogonal decompositions for such families were thoroughly
studied in a recent paper \cite{Kuznetsov3}. Relying mainly on the results of
{\it loc.cit.}, in Theorem \ref{th:flagexcept} we show that the derived categories of flag varieties ${\bf G}/{\bf
  B}$, where $\bf G$ is a semisimple Lie group of classical type, are generated by exceptional
collections. Note that Theorem \ref{th:flagexcept} (for classical
groups other than groups of type ${\bf A}_n$, which were treated in
the paper \cite{Kap} cited above) must have been known to experts for a long time;
however, there were no references to this fact.
\comment{
From the above theorems it follows that when looking for complete exceptional
collections on homogeneous spaces ${\bf G}/{\bf P}$, $\bf P$
being a parabolic subgroup of a semisimple Lie group $\bf G$, it is 
sufficient to restrict oneself to considering homogeneous
spaces that correspond to maximal parabolic subgroups of $\bf G$. 
}
In the final part of the paper we consider homogeneous spaces of the symplectic
groups of small rank, and prove in Theorem \ref{th:isotrgrass2-6}
that the derived category of homogeneous space ${\bf Sp}_6/{\bf
  P}$, where ${\bf P}$ corresponds to the middle node of the Dynkin
diagram of type ${\bf C}_3$, has a complete exceptional collection.
This gives, together with Theorem \ref{th:Samokh1} and the main result of
\cite{Sam1}, a complete description of the derived categories
of homogeneous spaces of the group ${\bf Sp}_6$ and of half a number of
homogeneous spaces of the group ${\bf Sp}_8$. Note also that the
homogeneous space ${\bf Sp}_6/{\bf P}$ in Theorem \ref{th:isotrgrass2-6} is
a smooth hyperplane section of the grassmannian ${\rm
  Gr}_{2,6}$. Recently, the derived categories of hyperplane sections of the
grassmannians ${\rm Gr}_{2,2n}$ for arbitrary $n$ 
were treated independently by A. Kuznetsov in \cite{Kuznetsov4}.\par

Finally, let us remark that a totally different approach to the above problems,
which is inspired by the work \cite{BMR}, is discussed in a
forthcoming paper \cite{Sam3}.

\section*{Notation}

All schemes are assumed to be smooth over an algebraically closed
field $\bf k$ of characteristic zero. 
For a proper scheme $X$ let $\Dd ^{b}(X)$ denote the bounded derived category of coherent
sheaves on $X$. For an object $\Ee$ of $\Dd ^{b}(X)$ we denote $\Ee
[n]$ the object $\Ee$ shifted by $n$. Given a morphism $f\colon
X\rightarrow Y$, we denote the derived functors of push--forward and pull--back
functors by $f_{\ast}$ and $f^{\ast}$, respectively. Similarly, the derived tensor product is denoted simply by
$\otimes$. For two objects $A$ and $B$ of $\Dd ^{b}(X)$ let ${\mathcal Hom}^{\bullet}_{X}(A,B)$ denote the object
${\mathcal R}{\mathcal Hom}_{X}(A,B) \in \Dd ^{b}(X)$ and 
let $\Hom ^{\bullet}_{X}(A,B)$ be the complex of $\bf k$-vector
spaces $\RHom _{\Dd ^{b}(X)}(A,B)$. Given a scheme $S$, points of 
$S$ always mean closed points. For a scheme $S$, and a vector bundle
$\Ee$ over $S$, the projectivization of $\Ee$ is denoted $\Pp
_{S}(\Ee)$. 

\section{Preliminaries}

We recall here some basic definitions. The material below is taken
from \cite{Bo}.

Throughout this section we fix a triangulated $\bf k$-linear category
${\sf D}$, equipped with a shift functor $[1]\colon {\sf D}\rightarrow {\sf D}$.  For two
objects $A, B \in {\sf D}$ let $\Hom ^{\bullet}_{\sf D}(A,B)$ be
the graded $\bf k$-vector space $\oplus _{i\in \mathbb Z}\Hom _{\sf
  D}(A,B[i])$. 
Let ${\sf A}\subset {\sf D}$ be a full triangulated subcategory,
that is a full subcategory of ${\sf D}$, which is closed
under shifts and taking cones.
\begin{definition}\label{def:orthogonalcat}
The right orthogonal ${\sf A}^{\perp}\subset \sf D$ is defined to be
the full subcategory
\begin{equation}
{\sf A}^{\perp} = \{B \in {\sf D}\colon \Hom _{\sf D}(A,B) = 0 \}
\end{equation}

\noindent for all $A \in {\sf A}$. The left orthogonal $^{\perp}{\sf
  A}$ is defined similarly. 

\end{definition}

\begin{definition}\label{def:admissible}
A full triangulated subcategory ${\sf A}$ of ${\sf D}$ is called
{\it right admissible} if the inclusion functor ${\sf A}\hookrightarrow {\sf
  D}$ has a right adjoint. Similarly, ${\sf A}$ is called {\it left
  admissible} if the inclusion functor has a left adjoint. Finally,
${\sf A}$ is {\it admissible} if it is both right and
left admissible.
\end{definition}

It is shown (\cite{Bo}, Lemma 3.1) that if a full triangulated category ${\sf
  A}\subset {\sf D}$ is right admissible then every object $X\in {\sf
  D}$ fits into a distinguished triangle
\begin{equation}
\dots \longrightarrow  Y\longrightarrow X\longrightarrow Z\longrightarrow Y[1]\rightarrow \dots
\end{equation}

\noindent with $Y\in {\sf A}$ and $Z\in {\sf A}^{\perp}$. One then
says that there is a semiorthogonal decomposition of ${\sf D}$ into
the subcategories $({\sf A}^{\perp}, \ {\sf A})$. More generally,
assume given a sequence of full triangulated subcategories ${\sf
  A}_1,\dots,{\sf A}_n \subset {\sf D}$. Denote $\langle {\sf
  A}_1,\dots,{\sf A}_n\rangle$ the triangulated subcategory of ${\sf
  D}$ generated by ${\sf A}_1,\dots,{\sf A}_n$.

\begin{definition}\label{def:semdecomposition}
A sequence $({\sf A}_1,\dots,{\sf A}_n)$ of admissible subcategories of
${\sf D}$ is called {\it semiorthogonal} if ${\sf
  A}_i\subset {\sf A}_j^{\perp}$ for $1\leq i < j\leq n$.
The sequence $({\sf A}_1,\dots,{\sf A}_n)$ is called a {\it semiorthogonal
  decomposition} of ${\sf D}$ if $\langle {\sf A}_1, \dots, {\sf A}_n
\rangle^{\perp} = 0$, that is ${\sf D} = \langle {\sf A}_1,\dots,{\sf A}_n\rangle$.
\end{definition}

\begin{definition}\label{def:exceptcollection}
An object $E \in \sf D$ is said to be exceptional if there is an
  isomorphism of graded $\bf k$-algebras 
\begin{equation}
\Hom _{\sf D}^{\bullet}(E,E) = {\bf k}.
\end{equation}

A collection of exceptional objects $(E_0,\dots,E_n)$ in $\sf D$ is called 
exceptional if for $1 \leq i < j \leq n$ one has
\begin{equation}
\Hom _{\sf D}^{\bullet}(E_j,E_i) = 0.
\end{equation}

\end{definition}

Denote $\langle E_0,\dots,E_n \rangle \subset {\sf D}$ the full
triangulated subcategory generated by the objects $E_0,\dots,E_n$. One 
proves (\cite{Bo}, Theorem 3.2) that such a category is admissible. 
The collection $(E_0,\dots,E_n)$ in $\sf D$ is said to be {\it complete} if 
$\langle E_0,\dots,E_n \rangle ^{\perp} = 0$, in other words ${\sf D}
= \langle E_0,\dots,E_n \rangle$.

\begin{definition}
Let $(E_0,E_1)$ be an exceptional pair in $\sf D$, i.e. an 
exceptional collection of two elements. The left mutation of $(E_0,E_1)$ is
a pair $(L_{E_0}E_1,E_0)$, where the object $L_{E_0}E_1$ is defined to
be a cone of the triangle
\begin{equation}
\dots \longrightarrow L_{E_0}E_1\longrightarrow \Hom _{\sf
  D}^{\bullet}(E_0,E_1)\otimes E_0\longrightarrow E_1\longrightarrow
  L_{E_0}E_1[1]\longrightarrow \dots \quad .
\end{equation}

\noindent The right mutation is a pair $(E_1,R_{E_1}E_0)$, where
  $R_{E_1}E_0$ is defined to be a cone of the triangle
\begin{equation}
\dots \longrightarrow R_{E_1}E_0[-1] \longrightarrow E_0\longrightarrow \Hom _{\sf
  D}^{\bullet}(E_0,E_1)^{\ast}\otimes E_1\longrightarrow
  R_{E_1}E_0\longrightarrow \dots \quad .
\end{equation}
\end{definition}

More generally, if $(E_0,\dots,E_n)$ is an exceptional collection of arbitrary length
in $\sf D$ then one can define left and right mutations of an object
$E\in {\sf D}$ through the category $\langle E_0,\dots,E_n
\rangle$. Denote $L_{\langle E_0,\dots,E_n
\rangle}E$ and $R_{\langle E_0,\dots,E_n \rangle}E$ left
and right mutations of $E$ through  $\langle E_0,\dots,E_n
\rangle$, respectively. One proves (\cite{Bo}, Proposition 2.1) that mutations of an exceptional collection are
exceptional collections.

Let $(E_0,\dots,E_n)$ be an exceptional collection in $\sf D$. One can
extend it to an infinite sequence of objects $(E_i)$ of $\sf D$, where
$i\in \mathbb Z$, defined inductively by putting 
\begin{equation}\label{eq:unwinding}
E_{i+n+1}\colon = R_{\langle E_{i+1},\dots,E_{i+n} \rangle}E_i, \qquad
E_{i-n-1}\colon = L_{\langle E_{i-n},\dots,E_{i-1} \rangle}E_i.
\end{equation}

Let $X$ be a variety of dimension $m$, and $\omega _X$ the
canonical invertible sheaf on $X$.
\begin{definition}
A sequence of objects $(E_i)$ of $\Dd ^{b}(X)$, where $i\in \mathbb Z$,
is called a {\it helix of period} $n$ if $E_i = E_{i+n}\otimes \omega
_{X}[m-n+1]$. An exceptional collection $(E_0,\dots,E_n)$ in $\Dd
^{b}(X)$ is called {\it a thread of the helix} if the infinite sequence
$(E_i)$, obtained from the collection $(E_0,\dots,E_n)$ as in (\ref{eq:unwinding}), is a helix of period $n+1$.
\end{definition}

For a Fano variety $X$ there is a criterion to establish whether a
given exceptional collection $(E_0,\dots,E_n)$ is complete, i.e. generates $\Dd
^{b}(X)$. 
\begin{theorem}[(\cite{Bo}, Theorem 4.1)]\label{th:Bondalth}
Let $X$ be a Fano variety, and $(E_0,\dots,E_n)$ be an exceptional
collection in $\Dd ^{b}(X)$. The following conditions are equivalent:
\begin{enumerate}

\item The collection $(E_0,\dots,E_n)$ generates $\Dd ^{b}(X)$.

\item The collection $(E_0,\dots,E_n)$ is a thread of the helix.

\end{enumerate}

\end{theorem}

\section{A theorem on semiorthogonal decompositions}\label{sec:genOrlovth}

The theorem below is a generalization of Orlov's theorem
on semiorthogonal decompositions in the derived categories of
projective bundles (\cite{Orlov}, Theorem 2.6).
\begin{theorem}\label{th:Samokh1}
Let $\pi : X\rightarrow S$ be a flat proper morphism 
  between two smooth schemes $X$ and $S$.
Let $\Ee _1,\dots,\Ee _n$ be a set of objects of $\Dd ^b(X)$ such that:

\begin{itemize}

\item [(i)] For any point $s\in S$ the restrictions
  $\Ee _{i}^s : = \Ee _{i}\otimes {\Oo}_{X_s}$, $X_s$ being the fiber
  $\pi ^{-1}(s)$, satisfy
  semiorthogonality condition, that is for $i > j$ one has $\Hom ^{\bullet}_{X_s}(\Ee _{i}^{s},\Ee
  _{j}^s) = 0$.

\item [(ii)] For $i = 1,\dots ,n$ and any point $s\in S$ the object
  $\Ee _{i}^s$ is an exceptional object of $\Dd ^{b}(X_s)$.

\item [(iii)] For any $s\in S$ the collection $(\Ee _{1}^s, \dots, \Ee _{n}^s)$ generates $\Dd ^{b}(X_s)$.

\end{itemize}

Put it simply, the above conditions say that for any $s\in S$ the
collection of objects $\Ee _{1}^s, \dots, \Ee _{n}^s$ is a complete exceptional
collection in $\Dd ^{b}(X_s)$.

Then $\Dd ^{b}(X)$ has a semiorthogonal decomposition:
\begin{equation}\label{eq:genOrlovth}
\Dd ^{b}(X) = \langle \pi ^{\ast}(\Dd ^b(S))\otimes \Ee
_1,\dots,\pi ^{\ast}(\Dd ^b(S))\otimes \Ee _n \rangle 
\end{equation}

\noindent Here $\pi ^{\ast}(\Dd ^{b}(S))\otimes \Ee _i$ denotes the full
triangulated subcategory in $\Dd ^{b}(X)$ generated by objects $\pi
^{\ast}(\Ff)\otimes \Ee _i, \ \Ff \in \Dd ^{b}(S)$. Each subcategory $\pi ^{\ast}(\Dd ^b(S))\otimes \Ee _i$ in the
decomposition (\ref{eq:genOrlovth}) is equivalent to $\Dd ^{b}(S)$. 

\end{theorem}

\begin{proof}
For $i = 1, \dots, n$ define functors $\Phi _i\colon \Dd
^{b}(S)\rightarrow \Dd ^{b}(X)$ by putting $\Phi _i(A) = \pi ^{\ast}(A)\otimes \Ee _i$,
where $A\in \Dd ^{b}(S)$. We check first that each $\Phi _i$ is a full
faithful embedding of $\Dd ^{b}(S)$ into $\Dd ^{b}(X)$.
\begin{lemma}
For $i = 1,\dots ,n$ the functors $\Phi _i\colon \Dd
^{b}(S)\rightarrow \Dd ^{b}(X)$ are full faithful embeddings.
\end{lemma}

\begin{proof}
Let $A,B\in \Dd ^{b}(S)$. For each $i = 1,\dots, n$ one has to show that
\begin{equation}\label{eq:functorsPhi}
\Hom _{X}^{\bullet}(\Phi _i(A),\Phi _i(B)) = \Hom _{S}^{\bullet}(A,B).
\end{equation}

\noindent By definition of the functors $\Phi _i$, the left-hand side of
(\ref{eq:functorsPhi}) is isomorphic to 
\begin{equation}\label{eq:Phiadjunction}
\Hom _{X}^{\bullet}(\pi ^{\ast}(A)\otimes \Ee _i,\pi ^{\ast}(B)\otimes
\Ee _i) = \Hom _{S}^{\bullet}(A,B\otimes \pi _{\ast}{\mathcal
  Hom}_{X}^{\bullet}(\Ee _i,\Ee _i)).
\end{equation}

\noindent The isomorphism (\ref{eq:Phiadjunction}) follows from adjunction of
the functors $\pi ^{\ast}$ and $\pi _{\ast}$, and projection
formula. Let us show that $\pi _{\ast}{\mathcal
  Hom}_{X}^{\bullet}(\Ee _i,\Ee _i)$ is quasiisomorphic to $\Oo _{S}$.
To this end, consider the diagram (cartesian square):
$$
\xymatrix @C5pc @R4pc {
 X_s \ar@{^{(}->}[r]^{{\tilde i}_s} \ar@<-0.1ex>[d]^{\pi _s} & X \ar@<-0.1ex>[d]^{\pi} \\
 s \ar@{^{(}->} @<-1pt>[r]^{i_s} &  S 
 }
$$

\noindent Since the morphism $\pi$ is flat, the base change holds 
(\cite{Kuznetsov1}, Corollary 2.23):
\begin{equation}\label{eq:sashabasechange}
{\pi _{s}}_{\ast}{\tilde i}_{s}^{\ast} = i_{s}^{\ast}\pi _{\ast}
\end{equation}

\noindent Hence, 
\begin{equation}\label{eq:basechange}
{\mathbb H}^{\ast}(X_s, {\tilde i}_{s}^{\ast}{\mathcal
  Hom}_{X}^{\bullet}(\Ee _i,\Ee _i)) = {\pi _{s}}_{\ast}{\tilde i}_{s}^{\ast}{\mathcal
  Hom}_{X}^{\bullet}(\Ee _i,\Ee _i) = i_{s}^{\ast}\pi _{\ast}{\mathcal
  Hom}_{X}^{\bullet}(\Ee _i,\Ee _i).
\end{equation}

\noindent By assumption (ii) of the theorem, for any $s\in S$ one has 
\begin{equation}\label{eq:basechangecontinued}
{\mathbb H}^{\ast}(X_s, {\tilde i}_{s}^{\ast}{\mathcal
  Hom}_{X}^{\bullet}(\Ee _i,\Ee _i)) = {\mathbb H}^{\ast}(X_s,{\mathcal
  Hom}_{X_s}^{\bullet}(\Ee _i^{s},\Ee _i^{s})) = {\bf k}.
\end{equation}

\noindent Thus, for any $s\in S$
\begin{equation}\label{eq:exceptobjrestr}
i_{s}^{\ast}\pi _{\ast}{\mathcal Hom}_{X}^{\bullet}(\Ee _i,\Ee _i) =
\Oo _s.
\end{equation}

\noindent Further, one has a canonical morphism $c\colon \Oo _{X}\stackrel{\rm
  c}\rightarrow {\mathcal Hom}_X^{\bullet}(\Ee _i,\Ee _i)$. 
Note that $ \Oo _{X} = \pi ^{\ast}\Oo _{S}$. Thus, by adjunction
one arrives at a distinguished triangle:
\begin{equation}\label{eq:pushforwcantriangle}
\dots \longrightarrow \Oo _{S}\stackrel{\tilde c}\longrightarrow \pi _{\ast}{\mathcal Hom}_X^{\bullet}(\Ee
_i,\Ee _i)\longrightarrow {\rm Cone}(\tilde c)\longrightarrow \Oo
_{S}[1]\longrightarrow \dots \ .
\end{equation}

\noindent Applying the functor $i_{s}^{\ast}$ to the triangle
(\ref{eq:pushforwcantriangle}) and taking into account
(\ref{eq:basechangecontinued}), and (\ref{eq:exceptobjrestr}), one gets 
\begin{equation}\label{eq:conerestriction}
i_{s}^{\ast}{\rm Cone}(\tilde c) = 0
\end{equation}

\noindent for any $s\in S$. Indeed, by flatness of $\pi$, the
morphism $i_{s}^{\ast}({\tilde c})\colon \Oo _{s}\rightarrow i_{s}^{\ast}
 \pi _{\ast}{\mathcal Hom}_X^{\bullet}(\Ee _i,\Ee _i)$, obtained from
 the triangle (\ref{eq:pushforwcantriangle}), is the isomorphism in
 (\ref{eq:basechangecontinued}). 
 Recall that the set $\Omega _{S} = \{\Oo _{s}: s\in S \}$, where $s$ ranges over the
set of all points in $S$, is a spanning class for $\Dd ^{b}(S)$
(\cite{Br1}, Example 2.2). It follows ({\it loc.cit.}) that an object $\Ff$ of $\Dd
^{b}(S)$ such that ${\rm Hom}_{S}^{\bullet}(\Ff, \Oo _{s}) = 0$ for
all $s\in S$ is quasiisomorphic to zero. For any $\Ff \in \Dd ^{b}(S)$ one has an isomorphism
\begin{equation}\label{eq:Ext-Torduality}
i_{s}^{\ast}{\Ff} = {\rm Hom}_{S}^{\bullet}(\Ff, \Oo _{s})^{\vee}.
\end{equation}

\noindent Taking $\Ff$ to be equal to ${\rm Cone}(\tilde
c)$ in (\ref{eq:Ext-Torduality}) and using (\ref{eq:conerestriction}),
one gets ${\rm Cone}(\tilde c) = 0$. Coming back to
(\ref{eq:pushforwcantriangle}) one sees that $\pi _{\ast}{\mathcal
  Hom}_{X}^{\bullet}(\Ee _i,\Ee _i) = \Oo _{S}$ , q.e.d.
\end{proof}

Put ${\sf D}_i\colon = \Phi _{i}(\Dd ^b(S))\subset \Dd ^{b}(X)$.
Let us check that each subcategory ${\sf D}_i$ is admissible in $\Dd
^{b}(X)$. 
\begin{lemma}
The subcategories ${\sf D}_1,\dots, {\sf D}_n\subset \Dd ^{b}(X)$ are admissible.
\end{lemma}

\begin{proof}
One has to make sure that the functors $\Phi _i\colon \Dd
^{b}(S)\hookrightarrow \Dd ^{b}(X)$ have left and right
adjoints. Indeed, for $A\in \Dd ^{b}(S)$ and $B\in \Dd ^{b}(X)$ one
has 
\begin{equation}\label{eq:rightadjforPhi}
\Hom _{X}^{\bullet}(\Phi _i(A),B) = \Hom _{X}^{\bullet}(\pi
^{\ast}(A)\otimes \Ee _i,B) = \Hom _{S}^{\bullet}(A,\pi
_{\ast}{\mathcal Hom}_{X}^{\bullet}(\Ee _i,B)),
\end{equation}

\noindent where the last isomorphism in (\ref{eq:rightadjforPhi})
follows from adjunction of the functors $\pi ^{\ast}$ and $\pi
_{\ast}$, and of the functors $\otimes$ and $\mathcal Hom$. Hence, a
right adjoint functor $\Phi _{i}^{!}$ to $\Phi _i$ is isomorphic to
$\pi _{\ast}{\mathcal Hom}_{X}^{\bullet}(\Ee _i,-)$. One checks
similarly that a left adjoint functor ${\Phi _i}_{\ast}$ to $\Phi _i$
is isomorphic to $\pi _{\ast}{\mathcal Hom}_{X}^{\bullet}(\Ee
_i\otimes \omega _{X/S}^{-1},-)[{\rm dim} X - {\rm dim} S]$, where
$\omega _{X/S}$ is the relative canonical invertible sheaf.

\end{proof}

We can prove now Theorem \ref{th:Samokh1}. Let us verify first the
semiorthogonality condition, i.e. that for any $k,l, \ k > l$ one has 
\begin{equation}\label{eq:verifsemcond}
\Hom _{X}^{\bullet}({\sf D}_k,{\sf D}_l) = \Hom _{X}^{\bullet}(\pi ^{\ast}(\Dd ^b(S))\otimes {\Ee}_k,\pi ^{\ast}(\Dd
(S))\otimes {\Ee}_l) = 0.
\end{equation}

\noindent Indeed, by adjunction of the functors $\pi ^{\ast}$ and $\pi
_{\ast}$ this is equivalent to
\begin{equation}\label{eq:vanishpushforw}
\pi _{\ast}{\mathcal Hom}_X^{\bullet}(\Ee _k,\Ee _l) = 0.
\end{equation}

\noindent Let $s\in S$ be a point. 
Using once again the base change along the diagram, one gets:
\begin{equation}\label{eq:basechange}
{\mathbb H}^{\ast}(X_s,{\tilde i}_{s}^{\ast}{\mathcal Hom}_X^{\bullet}(\Ee _k,\Ee
_l)) = {\pi _{s}}_{\ast}{\tilde i}_{s}^{\ast}{\mathcal Hom}_X^{\bullet}(\Ee _k,\Ee
_l) = i_{s}^{\ast}\pi _{\ast} {\mathcal Hom}_X^{\bullet}(\Ee _k,\Ee _l).
\end{equation}

\noindent By assumption (ii) of the theorem, the left-hand side of (\ref{eq:basechange})
is equal to zero. Therefore, for any $s\in S$ one has:
\begin{equation}\label{eq:restrvanish}
i_{s}^{\ast}\pi _{\ast}{\mathcal Hom}_X^{\bullet}(\Ee _k,\Ee _l) = 0.
\end{equation}

\noindent As above one gets $\pi _{\ast}{\mathcal Hom}_X^{\bullet}(\Ee _k,\Ee _l)
= 0$, the set $\Omega _{S}$ being a spanning class for $\Dd ^{b}(S)$. 
Hence, semiorthogonality is proven.\

To prove that the semiorthogonal collection of subcategories ${\sf
  D}_1,\dots,{\sf D}_n$ generates $\Dd ^b(X)$ it is
sufficient to show that $\langle {\sf D}_1,\dots,{\sf D}_n\rangle$ contains all the objects ${\Oo}_x, \ x\in
X$, the set $\Omega _{X} = \{\Oo _{x}: x\in X \}$ being a spanning class for $\Dd ^{b}(X)$. Every point
$x\in X$ lies in the fiber $X_s = \pi ^{-1}(s)$ for some $s\in S$. By
assumption (iii) of the theorem it follows that the sheaf $\Oo _{x}$ belongs to the
category generated by the objects $\Ee _{i}^{s}$:
\begin{equation}\label{eq:grattecielsheaves}
\Oo _{x}\in \langle \Ee _{1}^{s},\dots,\Ee _{n}^{s}\rangle .
\end{equation}

\noindent Recall that $\Ee _{i}^{s} = \Ee _{i}\otimes \Oo _{X_s}$ by
definition. Since $\pi ^{\ast}(\Oo _{s}) = \Oo _{X_s}$, the 
sheaf $\Oo _{x}$ belongs to $\langle {\sf D}_1,\dots,{\sf D}_n\rangle$
for any $x\in X$, therefore this semiorthogonal collection contains
the spanning class $\Omega _{X}$. This completes the proof.

\end{proof}

An immediate corollary to Theorem \ref{th:Samokh1} is a well-known fact:

\begin{corollary}
Let $X$ and $Y$ be two smooth projective varieties. Assume that
$(E_1,\dots ,E_n)$ and $(F_1,\dots ,F_m)$ are complete exceptional collections in ${\rm D}^{b}(X)$
and ${\rm D}^{b}(Y)$, respectively. Then $(E_1\boxtimes~F_1,\dots,E_{n}\boxtimes F_{m})$ is a complete exceptional collection in ${\rm
  D}^{b}(X\times Y)$.
\end{corollary}

\comment{
\begin{remark}
{\rm It follows {\it a posteriori} that under assumptions (ii) and
  (iii) of Theorem \ref{th:Samokh1} the morphism $\pi$ is, in fact,
  smooth. Indeed, \cite{BVdB}}
\end{remark}
}
\section{Complete exceptional collections on flag varieties}

In this section we show that the derived categories of flag
varieties of classical type are generated
by complete exceptional collections.
 
\begin{theorem}\label{th:flagexcept}
Let ${\bf G}$ be a semisimple algebraic group of classical type, $\bf B$ a
Borel subgroup, and ${\bf G}/{\bf B}$ the flag variety of $\bf G$. 
Then $\Dd ^{b}({\bf G}/{\bf B})$ is generated by an exceptional
collection. 
\end{theorem}        

\begin{proof}
Recall first that the case of homogeneous spaces of the groups of type ${\bf A}_n$ was
treated by M. Kapranov (\cite{Kap}, Theorem 3.10). We restrict ourselves, therefore, to the groups
of types ${\bf B}_n,{\bf C}_n$, and ${\bf D}_n$.
The proof is based on parabolic induction (e.g., \cite{LMS}), and we
refer the reader to {\it loc.cit.} for details of this construction.
Let us start with the symplectic group ${\bf Sp}_{2n}$, which
corresponds to the type ${\bf C}_n$, this case being the 
simplest. Fix a symplectic vector space $V$ of dimension $2n$, i.e. a space equipped with
a non-degenerate skew form $\omega \in \wedge ^{2}V^{\ast}$. 
The homogeneous space, which corresponds to the
leftmost node of the Dynkin diagram of ${\bf C}_n$, is isomorphic to
$\Pp (V) = \Pp ^{2n-1}$. Homogeneous spaces of the group ${\bf Sp}_{2n}$
are varieties of partial isotropic flags in $V$ with respect to $\omega$. The
variety of complete isotropic flags ${\rm IF}_{V}$, which is
isomorphic to ${\bf Sp}_{2n}/{\bf B}$, is equipped with a set of universal bundles
$\Uu _i$, where $i = 1, \dots, n$, fitting into a sequence:
\begin{equation}
0\subset {\Uu}_1\subset {\Uu}_2\subset \dots \subset {\Uu}_n =
{\Uu}_n^{\perp}\subset {\Uu}_{n-1}^{\perp}\subset \dots \subset {\Uu}_1^{\perp}\subset
V\otimes \Oo _{{\rm IF}_{V}}.
\end{equation}

\noindent Here ${\Uu}_i^{\perp}$ is a vector bundle that is
orthogonal to ${\Uu}_i$ with respect to $\omega$. The flag variety
${\rm IF}_{V}$ can be obtained as an iteration of projective
bundles over $\Pp (V)$. Indeed, for $m \leq n$ denote ${\rm
  IF}_{1,2,\dots,m}$ the partial flag variety of isotropic flags of type
$(1,2,\dots,m)$, the variety ${\rm IF}_{V}$ being
isomorphic to ${\rm IF}_{1,2,\dots ,n}$. 
By forgetting the top index, the variety ${\rm IF}_{1,2,\dots,m}$ becomes
the projectivization of symplectic vector bundle ${\mathcal N}_{m-1} = {\Uu}_{m-1}^{\perp}/{\Uu _{m-1}}$ over
${\rm IF}_{1,2,\dots,m-1}$. Arguing by induction on $m$, the base of induction being
$\Pp ^{2n-1} = {\rm IF}_{1}$, and applying Beilinson's theorem
(\cite{Be}), and Theorem \ref{th:Samokh1}, one obtains complete exceptional collections on
homogeneous spaces ${\rm IF}_{1,2,\dots,m}$ for $m\leq n$. 
For convenience of the reader, let us explicitly state the answer 
for the flag variety ${\rm IF}_{V}$ . Denote $\pi _{m} \colon
{\rm IF}_{1,2,\dots,m+1}\rightarrow {\rm
  IF}_{1,2,\dots,m}$ the projection, $\pi _{0}$ being the projection
to the point, and $p_{m}$ the composition
morphism $\pi _{m} \circ \dots \circ \pi _{n-1}\colon {\rm IF}_{1,2,\dots,n}\rightarrow {\rm
  IF}_{1,2,\dots,m}$. Let $\Oo _{{\mathcal N}_{m-1}}(-1)$ be the Grothendieck line bundle on $\Pp
_{{\rm IF}_{1,2,\dots,{m-1}}}({\mathcal N}_{m-1}) = {\rm
  IF}_{1,2,\dots,{m}}$. Then the collection of line bundles 
\begin{equation}
(\Oo _{{\mathcal N}_{n-1}}(j_{n-1})\otimes p_{n-1}^{\ast}\Oo _{{\mathcal
    N}_{n-2}}(j_{n-2})\otimes \dots \otimes p_{m}^{\ast}\Oo
    _{{\mathcal N}_{m-1}}(j_{m-1})\otimes \dots \otimes p_1^{\ast}\Oo _{\Pp ^{2n-1}}(j_0))
\end{equation}

\noindent is a complete exceptional collection in $\Dd
    ^{b}({\rm IF}_{V})$. Here the indices $j_k$ vary in the intervals ${\rm
    J}_k \colon = \{ -2n+2k < j_k \leq 0 \}$, the ordering on the set of
    indices $(j_{n-1},\dots, j_{0})$ being the product of linear
    orders on ${\rm J}_k$.\par

Consider now the case of orthogonal groups which correspond to the types ${\bf B}_n$ and ${\bf D}_n$.  
Fix an orthogonal vector space $V$ of dimension equal to either $2n+1$
or $2n$ for the groups of types ${\bf B}_n$ and ${\bf D}_n$,
respectively. The vector space $V$ is equipped with a non-degenerate
symmetric form $q\in {\sf S}^{2}V^{\ast}$. The homogeneous spaces, which correspond to the
leftmost nodes of the Dynkin diagrams of ${\bf B}_n$ and ${\bf
  D}_n$, are smooth quadrics in $\Pp (V)$ of dimensions
$2n-1$ and $2n-2$, respectively. The derived category of a quadric
was studied in \cite{Kap}. Recall its description:

\begin{theorem}[(\cite{Kap}, Theorem 4.10)]\label{th:Kapthquadrics}
Let ${\rm Q}\subset \Pp ^{n+1}$ be a smooth quadric of dimension
$n$. Then $\Dd ^{b}({\rm Q})$ is generated by an exceptional
collection
\begin{equation}
\Sigma (-n), {\Oo}_{Q}(-n+1),\dots, {\Oo}_{Q}(-1),{\Oo}_{Q}  
\end{equation}

if $n$ is odd, and 
\begin{equation}
\Sigma _{+}(-n), \Sigma _{-}(-n), {\Oo}_{Q}(-n+1),\dots, {\Oo}_{Q}(-1),{\Oo}_{Q} 
\end{equation}

if $n$ is even.

\end{theorem}

\noindent Here $\Sigma _{-},\Sigma _{+}$ are spinor bundles (they are isomorphic
to each other if $n$ is odd and the isomorphism class is denoted by $\Sigma$).\par

The homogeneous spaces of orthogonal group of either type 
are varieties of partial isotropic flags in $V$ with respect to $q$.
We restrict ourselves to the case of the orthogonal group of type ${\bf
  B}_n$, the case of ${\bf D}_n$ being similar. 
The variety of complete isotropic flags ${\rm OF}_{V}$, which is
isomorphic to ${\bf SO}_{2n+1}/{\bf B}$, is equipped with a set of universal bundles
$\Ww _i$, where $i = 1, \dots, n$, fitting into a sequence:
\begin{equation}
0\subset {\Ww}_1\subset {\Ww}_2\subset \dots \subset {\Ww}_n \subset
{\Ww}_n^{\perp}\subset {\Ww}_{n-1}^{\perp}\subset \dots \subset {\Ww}_1^{\perp}\subset
V\otimes \Oo _{{\rm IF}_{V}}.
\end{equation}

\noindent  Here ${\Ww}_i^{\perp}$ is a vector bundle that is
orthogonal to ${\Ww}_i$ with respect to $q$.
An argument, which is analogous to that outlined for the
group ${\bf Sp}_{2n}$, shows that the flag variety ${\rm OF}_{V}$ can be
obtained as a successive iteration of smooth quadric fibrations over a smooth
quadric ${\rm Q}_{2n-1}\subset \Pp (V)$. Unfortunately, Theorem \ref{th:Samokh1}
cannot be applied directly to such
fibrations, the reason being that spinor bundles are not globally
defined in general. One could prove still that in our particular case of quadric
fibrations spinor bundles do exist globally and then use Theorem
\ref{th:Samokh1}. It makes little sense, however, since a general theorem on semiorthogonal
decompositions for flat families of quadrics is now available
(\cite{Kuznetsov3}). Recall briefly the setup of {\it loc.cit}. Let $S$ be a smooth scheme
(the base of a family), and let $\Ee$ be a vector bundle of rank $n$
on $S$. Denote by $\Oo _{\Ee}(1)$ the Grothendieck invertible sheaf on
$\Pp _{S}(\Ee)$. Fix a line bundle $\Ll$ on $S$, and let $\varphi \colon {\sf
  S}^{2}\Ee \rightarrow \Ll$ be a symmetric bilinear form on $\Ee$
that makes $\Ee$ into an orthogonal vector bundle. 
The form $\varphi$ defines a section of the line bundle $\Oo
_{\Ee}(2)\otimes \pi ^{\ast}\Ll$ on $\Pp _{S}(\Ee)$. Let $X$ be the zero locus of this
section, $i$ the embedding $X \hookrightarrow \Pp
_{S}(\Ee)$, and $p$ the composition morphism $\pi \circ i \colon
X\rightarrow S$. Finally, denote by $\Oo _{X/S}(1)$ the 
line bundle $i^{\ast}\Oo _{\Ee}(1)$. The projection $p$ makes $X$ into a quadric
fibration over $S$. One checks, as in Theorem \ref{th:Samokh1}, that for $k
= 1,\dots, n - 2$ the collection of admissible subcategories $p^{\ast}\Dd
^{b}(S)\otimes \Oo _{X/S}(k)\subset \Dd ^{b}(X)$ is semiorthogonal.
The right orthogonal in $\Dd ^{b}(X)$ to this collection is described by 
\begin{theorem}[(\cite{Kuznetsov3}, Theorem 4.2)]\label{th:SashaKquadrics}
There exists a semiorthogonal decomposition in $\Dd ^{b}(X)$:
\begin{eqnarray}\label{eq:Kuzdecomforquad}
& \Dd ^{b}(X) = \langle \Dd ^{b}(S,{\mathcal B}_{0}),
p^{\ast}\Dd^{b}(S)\otimes \Oo _{X/S}(1), p^{\ast}\Dd ^{b}(S)\otimes
\Oo _{X/S}(2), \dots, & \nonumber\\
& p^{\ast}\Dd ^{b}(S)\otimes \Oo _{X/S}(n-2)\rangle
\end{eqnarray}

\noindent where $\Dd ^{b}(S,{\mathcal B}_{0})$ is the derived
category of coherent sheaves of ${\mathcal B}_0$-modules over $S$.

\end{theorem}

Here ${\mathcal B}_0$ is the sheaf of even parts of 
Clifford algebras on $S$ ({\it loc.cit.}, Section 3.1).\par

\noindent To complete our proof, we argue again by induction, the
base of induction being the quadric ${\rm Q}_{2n-1}$. 
Indeed, for $m \leq n$ denote ${\rm
  OF}_{1,2,\dots,m}$ the partial flag variety of isotropic flags of type
$(1,2,\dots,m)$, the flag variety ${\rm OF}_{1,2,\dots,n}$ being
isomorphic to ${\rm OF}_{V}$. Denote $\pi _{m} \colon
{\rm OF}_{1,2,\dots,m+1}\rightarrow {\rm
  OF}_{1,2,\dots,m}$ the projection, $\pi _{0}\colon {\rm
  Q}_{2n-1}\rightarrow {\rm pt}$ being the projection
to the point, and let $p_{m}$ be the composition
morphism $\pi _{m} \circ \dots \circ \pi _{n-1}\colon {\rm OF}_{1,2,\dots,n}\rightarrow {\rm
  OF}_{1,2,\dots,m}$. The variety ${\rm OF}_{1,2,\dots,m}$ becomes then a smooth quadric fibration 
over ${\rm OF}_{1,2,\dots,m-1}$. If ${\mathcal M}_{m-1} =
{\Ww}_{m-1}^{\perp}/{\Ww _{m-1}}$ 
then there is an embedding $i_{m}\colon {\rm OF}_{1,2,\dots,m}\subset
\Pp _{{\rm OF}_{1,2,\dots,m-1}}({\mathcal M}_{m-1})$.

Theorem \ref{th:SashaKquadrics} furnishes a semiorthogonal decomposition
of $\Dd ^{b}({\rm OF}_{1,2,\dots,m})$. Denote ${\mathcal B}_0^{m}$ the
corresponding sheaf of even parts of Clifford algebras on ${\rm OF}_{1,2,\dots,m-1}$.
The structure of the sheaf ${\mathcal
  B}_{0}^{m}$ gets more simple in this case, the
morphism $\pi _{m-1}$ being smooth. Indeed, Proposition 3.15 of 
\cite{Kuznetsov3} ensures that the sheaf ${\mathcal B}_0^{m}$ is a sheaf of
Azumaya algebras over ${\rm OF}_{1,2,\dots,m-1}$. Moreover, ${\mathcal B}_0^{m}$
splits as a sheaf of Azumaya algebras, as shows the following lemma.

\begin{lemma}\label{prop:Brauergroup}
Let $X$ be a variety admitting a cellular algebraic
decomposition. Then the  cohomological Brauer group ${\rm Br}^{'}(X) =
{\rm H}^{2}(X,\Oo _{X}^{\ast})_{\rm tors}$ is trivial.
\end{lemma}

\begin{proof}
It is known that for such a variety the cycle map ${\rm cl}_{X}$ provides an
  isomorphism between Chow groups ${\rm A}_{\ast}(X)$ and 
  homology groups ${\rm H}_{\ast}(X)$ (\cite{Fu}, Example
  19.1.11). Thus, all cycles on $X$ are algebraic and 
  odd-dimensional integral cohomology groups of $X$ are trivial. 
Considering exponential short exact sequence, one gets an 
exact sequence for ${\rm Br}^{'}(X)$ (\cite{Sch}, Proposition 1.1):
\begin{equation}
0\rightarrow A\otimes {\mathbb Q}/{\mathbb Z} \rightarrow {\rm
  Br}^{'}(X)\rightarrow T\rightarrow 0
\end{equation}

\noindent where $T$ is the torsion part of ${\rm H}^{3}(X,\mathbb Z)$, and $A$ is the
lattice of transcendental cycles, which is isomorphic to 
${\rm H}^{2}(X,{\mathbb Z})/{\rm NS}(X)$, the factor group of ${\rm
  H}^{2}(X,{\mathbb Z})$ by the Neron--Severi group. 
Given the above remarks, we get ${\rm Br}^{'}(X) = 0$.
\end{proof}

By \cite{Gro}, the Brauer group ${\rm Br}(X)$, which parametrizes
equivalence classes of sheaves of Azumaya algebras on $X$, naturally
injects into ${\rm Br}^{'}(X)$. Homogeneous spaces of semisimple Lie
groups have algebraic cellular
decompositions by virtue of the Bruhat decomposition. By Lemma \ref{prop:Brauergroup}, the
group ${\rm Br}^{'}({\rm OF}_{1,2,\dots,m-1})$ is trivial, thus the
corresponding sheaf of Azumaya
algebras ${\mathcal B}_0^{m}$ over ${\rm OF}_{1,2,\dots,m-1}$ splits,
i.e. is isomorphic to ${\mathcal
  End}(\Sigma _m)$ for a vector bundle $\Sigma _m$. If $m=1$, i.e. one has
the quadric ${\rm Q}_{2n-1}$, then Proposition 4.4 of \cite{Kap}
asserts that a splitting bundle for ${\mathcal B}_0^{1}$ is isomorphic to the
spinor bundle $\Sigma = \colon \Sigma _1$. In general, we have relative spinor bundles
$\Sigma _{m}$ on ${\rm OF}_{1,2,\dots,m-1}$, splitting the sheaves of
Azumaya algebras 
${\mathcal B}_0^{m}$. Finally, the category $\Dd ^{b}(S,{\mathcal
  A})$ for a split Azumaya algebra $\mathcal A$ is equivalent to $\Dd ^{b}(S)$,
since abelian categories ${\rm Coh}(S)$ and ${\rm Coh}({\mathcal
  A})$ are equivalent in this case 
(the equivalence being given by the map $\Ff \rightarrow {\sf E}\otimes
\Ff$, where $\Ff \in {\rm Coh}(S)$ and $\sf E$ is a splitting bundle
for $\mathcal A$). Moreover, an embedding of $\Dd ^{b}({\rm
  OF}_{1,2,\dots,m-1},{\mathcal B}_{0}^{m})\simeq \Dd ^{b}({\rm
  OF}_{1,2,\dots,m-1})$ into $\Dd ^{b}({\rm OF}_{1,2,\dots,m})$ can be given
by the functor $\pi _{m-1}^{\ast}(-)\otimes \Sigma _{m}$
(\cite{Kuznetsov3}, Proposition 4.9).
Arguing by induction on $m$ and applying Theorems \ref{th:Kapthquadrics}, \ref{th:SashaKquadrics}, and
the above arguments, one obtains, therefore, complete exceptional collections on
homogeneous spaces ${\rm OF}_{1,2,\dots,m}, m\leq n$. 
For convenience of the reader, let us explicitly state the answer 
for the flag variety. Let $\Oo _{{\mathcal M}_{m-1}}(-1)$ be the Grothendieck line bundle on $\Pp
_{{\rm OF}_{1,2,\dots,{m-1}}}({\mathcal M}_{m-1})$ and let $\Oo _{{\sf
    M}_{m-1}}(-1)$ be its restriction to ${\rm OF}_{1,2,\dots,{m}}$. Then the collection of vector bundles 
\begin{eqnarray}\label{eq:exccollectortflags}
& ( \Sigma _{n}\otimes p^{\ast}_{n-1}\Sigma _{n-1}\otimes \dots \otimes
p_{m}^{\ast}\Sigma _{m}\otimes \dots \otimes p_{1}^{\ast}\Sigma _{1}, \dots ,& \nonumber \\
&   p^{\ast}_{n-1}\Sigma _{n-1}\otimes \dots \otimes
p_{m}^{\ast}\Sigma _{m}\otimes p_{m}^{\ast}\Oo _{{\sf M}_{m}}(j_{m})\otimes \dots
\otimes p_{1}^{\ast}\Oo _{{\rm Q}_{2n-1}}(j_{0}), \dots & \nonumber \\
& \Oo _{{\sf M}_{n-1}}(k_{n-1})\otimes \dots \otimes p_{m}^{\ast}\Oo _{{\sf M}_{m}}(k_{m})\otimes \dots
\otimes p_{1}^{\ast}\Oo _{{\rm Q}_{2n-1}}(k_{0}), \dots , \Oo _{{\rm OF}_V} ) , & 
\end{eqnarray}

\noindent where $0 \leq j_m, k_m \leq 2n - 2m - 2$ and $0\leq m \leq n
-1$, is a complete exceptional
collection in $\Dd ^{b}({\rm OF}_{V})$. Note that $\Sigma _{n}$ is
isomorphic to $\Oo _{{\sf M}_{n-1}}(-1)$. The leftmost bundle in
(\ref{eq:exccollectortflags}) is the tensor product of 
all spinor bundles $\Sigma _i$ pulled back to ${\rm OF}_{V}$ for
$i=1,\dots ,n$, and the rightmost bundle is the structural sheaf $\Oo
_{{\rm OF}_V}$. Terms in the middle line of
(\ref{eq:exccollectortflags}) are pull-backs to $\Oo _{{\rm OF}_V}$ of
tensor products of the spinor bundles and line bundles. Finally,
terms in the bottom line are pull-backs to $\Oo _{{\rm OF}_V}$ of
tensor products of line bundles. The ordering of the collection (\ref{eq:exccollectortflags}) is a 
natural one obtained from the construction.
\comment{
The terms of the collection (\ref{eq:exccollectortflags}) are ordered as
follows: a bundle ${\sf E}_1$ is to left 
}
\end{proof}

\begin{remark}
{\rm  Note that the collections constructed on the flag varieties for
  orthogonal groups do not consist entirely of line bundles. Conjecturally, there exists a
complete exceptional collection of line bundles on the flag variety of
a semisimple Lie group. One would like to have a uniform construction of
such collections on flag varieties.}
\end{remark}

\section{Homogeneous spaces of the symplectic groups of small
  ranks}\label{sec:symplgroups}

In this section we study homogeneous spaces of Lie groups of type
${\bf C}_n$ for $n\leq 4$. The main result here
is Theorem \ref{th:isotrgrass2-6} where we find a complete exceptional
collection in the derived category of a homogeneous space of the group
${\bf Sp}_6$. Theorems \ref{th:isotrgrass2-6} and \ref{th:Samokh1},
and the main result of \cite{Sam1}, yield complete exceptional
collections in the derived categories of all homogeneous spaces of the group
${\bf Sp}_6$. These results also make possible to find complete exceptional
collections on some homogeneous spaces of ${\bf Sp}_8$. \par 

Consider the group ${\bf Sp}_6$. Its Dynkin diagram is of type ${\bf
C}_3$ and there are a total of seven homogeneous spaces of ${\bf 
Sp}_6$. Fix a symplectic vector space $V$ of dimension 6, and let $\omega$ be a non-degenerate skew
form on $V$, $\omega \in \wedge ^{2}V^{\ast}$. 
Maximal parabolic subgroups of ${\bf Sp}_6$ correspond to the nodes of the diagram ${\bf C}_3$. The homogeneous space,
corresponding to the leftmost node is $\Pp ^5$, and the homogeneous
space, corresponding to the rightmost node is ${\rm LG}_{3,6}$, the
grassmannian of lagrangian subspaces in $V$. The category $\Dd 
^{b}({\rm LG}_{3,6})$ was described in \cite{Sam1}. Below we treat
the case of homogeneous space, corresponding to the middle node of
${\bf C}_3$.

This is the variety ${\rm IGr}_{2,6}$ of two-dimensional isotropic subspaces in $V$ with
respect to $\omega$. It is easy to see that ${\rm IGr}_{2,6}$ is a smooth hyperplane section of the
grassmannian ${\rm Gr}_{2,6}$. Indeed, ${\rm Gr}_{2,6}$ is embedded in $\Pp (\wedge ^{2}V^{\ast})$ via Pl\"ucker
embedding. The form $\omega$ defines a section $s_{\omega}$ of the 
corresponding line bundle $\Oo (1)$ on ${\rm Gr}_{2,6}$. The zero locus of
$s_{\omega}$, that is a hyperplane section of ${\rm Gr}_{2,6}$ with respect to the above
embedding, is the locus of two-dimensional subspaces in $V$ such that $\omega$
vanishes on these subspaces, that is ${\rm IGr}_{2,6}$. The dimension of
${\rm IGr}_{2,6}$ is therefore equal to 7. 
Let ${\rm IF}_{1,2,6}$ be the space of pairs $(l\subset U)$, where $l$
and $U$ are a one-dimensional and an isotropic two-dimensional vector subspaces in
$V$, respectively. One has a diagram (with obvious projections):

\vspace*{20pt}
$$
\xymatrix @C5pc @R4pc {
 & {\rm IF}_{1,2,6} \ar[dl]_{p} \ar@<-0.1ex>[dr]^{q} &   \\
 \Pp ^{5}           &                    &           {\rm IGr}_{2,6}
 }
$$
\vspace*{20pt}

The projections $p$ and $q$ make ${\rm IF}_{1,2,6}$ into a $\Pp
^3$-bundle over $\Pp ^{5}$ and a $\Pp ^1$-bundle over ${\rm
  IGr}_{2,6}$, respectively. More exactly, the $\Pp ^3$-bundle over $\Pp ^{5}$ is the
projectivization of the null-correlation bundle ${\sf N}_{\omega}$ of rank 4 that
corresponds to the form $\omega$ (\cite{OSS}, $\S 4$). Recall that ${\sf
  N}_{\omega}$ fits into a short exact sequence of vector bundles on
$\Pp ^{5}$:
\begin{equation}
0\longrightarrow {\sf N}_{\omega}\longrightarrow {\mathcal T}_{\Pp
  ^5}(-1)\stackrel{p_{\omega}}\longrightarrow \Oo _{\Pp ^5}(1)\longrightarrow 0 ,
\end{equation}

\noindent where ${\mathcal T}_{\Pp ^5}$ is the tangent bundle on $\Pp
  ^5$, and the map $p_{\omega}\in \Hom _{\Pp ^5}({\mathcal
  T}_{\Pp^5}(-1),\Oo _{\Pp ^5}(1)) = \bigwedge ^{2}V^{\ast}$
  corresponds to the form $\omega$.
The $\Pp ^{1}$-bundle over ${\rm IGr}_{2,6}$ is the projectivization of 
the universal bundle $\Uu$ of rank 2. From the above diagram one
easily finds the rank of the Grothendieck group ${\rm
  K}^{0}({\rm IGr}_{2,6})$ to be equal to 12. Given such a diagram, 
it is natural to suppose that the direct images under the
morphism $q$ of appropriate line bundles on ${\rm IF}_{1,2,6}$ may give a
complete exceptional set on ${\rm IGr}_{2,6}$ (cf. \cite{Sam1}). This
method works in this case indeed. 

\begin{theorem}\label{th:isotrgrass2-6}
$\Dd ^b({\rm IGr}_{2,6})$ is generated by the exceptional
  collection:
\begin{eqnarray}
& {\sf K} = \langle \Uu(-4), \ \Oo _{\rm IGr_{2,6}}(-4), \ {\Ss}^2\Uu(-3),  \ \Uu(-3), \ \Oo
 _{\rm IGr_{2,6}}(-3), \ {\Ss}^2\Uu(-2),\Uu(-2), & \nonumber \\
& \Oo _{\rm IGr_{2,6}}(-2), \Uu(-1),\Oo _{\rm IGr_{2,6}}(-1), \Uu,\Oo
 _{\rm IGr_{2,6}}\rangle .
\end{eqnarray}
\end{theorem}

\begin{proof}
One immediately verifies, using the Borel--Weil--Bott theorem,
that the collection $\sf K$ is exceptional. We need to prove that $\sf K$
is complete, i.e. generates the derived category.\par

Let $\Oo_{\Uu}(-1)$ and ${\Oo}_{{\sf
    N}_{\omega}}(-1)$ be the Grothendieck invertible sheaves
on $\Pp _{\rm IGr_{2,6}}(\Uu )$ and $\Pp _{\Pp ^5}({\sf
    N}_{\omega})$, respectively. Denote $\omega _{\rm IF_{1,2,6}}$ the
    canonical line bundle on ${\rm IF_{1,2,6}}$. There are following
    relations between line bundles on ${\rm IF_{1,2,6}}$: 

\begin{itemize}

\item[(i)] $p^{\ast}{\Oo}_{\Pp ^5}(-1) = {\Oo}_{\Uu}(-1)$

\item[(ii)] $p^{\ast}{\Oo}_{\Pp ^5}(-1)\otimes {\Oo}_{{\sf N}_{\omega}}(-1) =
  q^{\ast}{\Oo}_{\rm IGr_{2,6}}(-1)$

\item[(iii)] $\omega _{\rm IF_{1,2,6}} = p^{\ast}{\Oo}_{\Pp ^{5}}(-2)\otimes
  q^{\ast}{\Oo}_{\rm IGr_{2,6}}(-4) = p^{\ast}{\Oo}_{\Pp
  ^{5}}(-6)\otimes \Oo _{{\sf N}_{\omega}}(-4)$

\end{itemize}

\noindent Let ${\Ll}_{i,j}$ be the sheaf $p^{\ast}{\Oo
}_{\Pp ^{5}}(-i)\otimes {\Oo}_{\sf N _{\omega}}(-j)$. From the above relations one
finds:


\begin{equation}\label{fig:fig11}
{q}_{\ast}({\Ll}_{i,j}) = \left\{
                          \begin{aligned}
                          {\Ss}^{j-i}{\Uu}(-i), \quad j\geq i \\
                          0, \quad j - i = -1 \\
                          {\Ss}^{i-j-2}{\Uu}(-j-1)[-1], \quad j - i \leq -2\\
                          \end{aligned}
                          \right. \quad .
\end{equation}

\vspace*{0.2cm}

\noindent For $i=-3,\dots,0$ denote ${\sf D}_i$ the full subcategory $p^{\ast}\Dd
^{b}(\Pp ^{5})\otimes \Oo _{\sf N_{\omega}}(i)\subset \Dd ^{b}({\rm
  IF}_{1,2,6})$. There is a semiorthogonal decomposition in $\Dd ^{b}({\rm
  IF}_{1,2,6})$ with respect to the projection $p$:
\begin{eqnarray}\label{eq:apprsemdec}
& \Dd ^{b}({\rm IF}_{1,2,6}) = \langle {\sf D}_{-3}, {\sf D}_{-2},
{\sf D}_{-1}, {\sf D}_{0} \rangle .
\end{eqnarray}

To prove the theorem we need to pick up a convenient complete
exceptional collection in $\Dd ^{b}({\rm IF}_{1,2,6})$. To this end,
choose the following exceptional collections in ${\sf D}_i$ for $i=-3,\dots,0$:
\begin{eqnarray}
& {\sf D}_{-3} = \langle \Oo _{\Pp ^{5}}(-7), \Oo _{\Pp ^{5}}(-6),
  \dots, \Oo _{\Pp ^{5}}(-3), \Oo _{\Pp ^{5}}(-2)\rangle \otimes \Oo _{\sf
  N_{\omega}}(-3) & \nonumber \\
& {\sf D}_{-2} = \langle \Oo _{\Pp ^{5}}(-6), \Oo _{\Pp ^{5}}(-5),
  \dots, \Oo _{\Pp ^{5}}(-2), \Oo _{\Pp ^{5}}(-1)\rangle \otimes \Oo _{\sf
  N_{\omega}}(-2) & \nonumber \\
& {\sf D}_{-1} = \langle \Oo _{\Pp ^{5}}(-5), \Oo _{\Pp ^{5}}(-4),
  \dots, \Oo _{\Pp ^{5}}(-1), \Oo _{\Pp ^{5}}\rangle \otimes \Oo _{\sf
  N_{\omega}}(-1) & \nonumber \\
& {\sf D}_{0} = \langle \Oo _{\Pp ^{5}}(-4), \Oo _{\Pp ^{5}}(-3),
  \dots, \Oo _{\Pp ^{5}}, \Oo _{\Pp ^{5}}(1)\rangle \otimes \Oo _{\sf
  N_{\omega}} .
\end{eqnarray}

\noindent Denote ${\sf L}$ the complete exceptional collection in $\Dd
  ^{b}({\rm IF}_{1,2,6})$ thus obtained. Let $\Mm$ be the set of objects of $\Dd ^{b}({\rm IGr}_{2,6})$ obtained by
applying the functor $q_{\ast}$ to the collection ${\sf
  L}$. Write this symbolically as 
\begin{equation}\label{eq:surjectivity}
q_{\ast} \ {\sf L} = \Mm .
\end{equation}

\noindent From (\ref{fig:fig11}) one sees that up to a shift all the sheaves 
of the collection $\sf K$
are contained in $\Mm$. If the sets $\sf K$ and $\Mm$ consisted of the same elements
this would prove the theorem. This is not the case,
however, since two objects of $\Mm$ do not belong to $\sf K$:
\begin{equation}\label{eq:notincollection}
 {q}_{\ast}({\Ll}_{7,3}) = {\Ss}^{2}{\Uu}(-4)[-1] \qquad \mbox{and} \qquad
{q}_{\ast}({\Ll}_{-1,0}) = {\Uu}^{\ast}.
\end{equation}

\noindent The invertible sheaves ${\Ll}_{7,3}$ and ${\Ll}_{-1,0}$ are situated at the opposite ends of 
  the collection ${\sf L}$. We can mutate $\Ll _{7,3}$ to
  the right through ${\sf L}$. By Theorem \ref{th:Bondalth}, this
  mutation is given, up to a shift, by  
  twisting the bundle $\Ll _{7,3}$ by $\omega _{\rm IF_{1,2,6}}^{-1} = \Ll _{-6,-4}$. 
We obtain, therefore, the line bundle $\Ll _{1,-1}$. Such a mutation 
gives rise to a new thread of the helix generated by $\sf L$. With
  respect to this thread the category $\Dd ^{b}({\rm IF}_{1,2,6})$
  has a semiorthogonal decomposition:
\begin{equation}
\Dd ^{b}({\rm IF}_{1,2,6}) = (\langle \Ll _{-1,0}, \Ll _{1,-1} \rangle ^{\perp},
\langle \Ll _{-1,0}, \Ll _{1,-1} \rangle ).
\end{equation}

\noindent Assume that the collection $\sf K$ is not complete, i.e. there
is a non-zero left orthogonal $^{\perp}{\sf K} \neq 0$. Let $\Ff$ be a
non-zero object of $^{\perp}{\sf K}$. Then $q^{\ast}\Ff$ 
belongs to the full subcategory $\langle \Ll _{-1,0}, \Ll _{1,-1}
\rangle$. Indeed, for any exceptional generator $\mathcal M$ of the
category $\langle \Ll _{-1,0}, \Ll _{1,-1} \rangle ^{\perp}$, one
has:
\begin{equation}
\Hom ^{\bullet}_{\rm IF_{1,2,6}}(q^{\ast}\Ff, {\mathcal M}) = \Hom
^{\bullet}_{\rm IGr_{2,6}}(\Ff, q_{\ast}{\mathcal M}) =  0,
\end{equation}

\noindent since $q_{\ast}{\mathcal M}$ belongs to $\sf K$ by the very
construction. Hence, there is a distinguished triangle:
\begin{equation}\label{eq:trian1}
\dots \longrightarrow V_{1,-1}^{\cdot}\otimes {\Ll _{1,-1}} \longrightarrow
  q^{\ast}\Ff \longrightarrow V_{-1,0}^{\cdot}\otimes {\Ll _{-1,0}} \longrightarrow
  \dots \ .
\end{equation}

\noindent Applying the functor $q_{\ast}$ to the triangle (\ref{eq:trian1}), one gets:
\begin{equation}\label{eq:trian2}
\dots \longrightarrow V_{1,-1}^{\cdot}\otimes \Oo _{\rm IGr_{2,6}}[-1]
\longrightarrow \Ff \longrightarrow V_{-1,0}^{\cdot}\otimes \Uu
^{\ast} \longrightarrow \dots \ .
\end{equation}

\noindent Since $\Hom ^{\bullet}_{\rm IGr_{2,6}}(\Uu ^{\ast},\Oo _{\rm
  IGr_{2,6}}) =
0$, the triangle (\ref{eq:trian2}) splits and one has an isomorphism:
\begin{equation}\label{eq:split}
\Ff = V_{-1,0}^{\cdot}\otimes \Uu ^{\ast}\oplus
V_{1,-1}^{\cdot}\otimes \Oo _{\rm IGr_{2,6}}[-1].
\end{equation}

\noindent The \ graded \ vector \ space
\ $V_{1,-1}$ \ is \ forced \ to \ be \ zero, \ since \ $\Ff$ $\in$ 
$^{\perp}{\sf K}$ \ and
$\Hom ^{\bullet}_{\rm IGr_{2,6}}(\Uu ^{\ast},\Oo _{\rm IGr_{2,6}}) =
0$. Therefore, $\Ff$ is isomorphic to
$V_{-1,0}^{\cdot}\otimes \Uu _{2}^{\ast}$. Note that $\Hom
^{\bullet}_{\rm IGr_{2,6}}(\Uu ^{\ast},\Uu (-4)) = {\bf
  k}[-7]$. Taking $\Hom$-s between $\Ff$ and various shifts of the bundle $\Uu (-4)$, and using the orthogonality
condition, one sees that the graded vector space $V_{-1,0}^{\cdot}$ is
zero as well. Hence, $\Ff$ is zero. This finishes the proof.
\end{proof}

It follows easily now that the derived category of any homogeneous
space for ${\bf Sp}_6$ is generated by an exceptional
collection. Indeed, such a homogeneous space is a variety of partial
isotropic flags in $V$, and is obtained by giving an ordered subset $S$ of
the ordered set $(1,2,3)$. Denote this variety $\rm IF_{S}$. The derived
categories of ${\rm IF}_{i}$ for $i = 1,2,3$, by \cite{Be}, \cite{Sam1}, and Theorem
\ref{th:isotrgrass2-6} are generated by exceptional collections. Any
other homogeneous space ${\rm IF}_{S}$ is mapped onto one of ${\rm
  IF}_{i}$ for $i = 1,2,3$, and it is easy to see that ${\rm IF}_{S}$ is
either a $\Pp ^{1}$-bundle over ${\rm IF}_{i}$ for some $i$, or a $\Pp
^{2}$-bundle over ${\rm IF}_{3}$, or a 
two-step iteration of $\Pp ^{k}$-bundles for $k = 1$ and $2$ (the last case corresponds to the flag
variety for ${\bf Sp}_6$). Therefore, applying either Orlov's theorem
(\cite{Orlov}) or Theorem \ref{th:Samokh1}, and the above arguments
one gets complete exceptional collections in the derived categories of
homogeneous spaces of the group ${\bf Sp}_6$.
\par 

Arguing as above, we also obtain that the derived categories of some
homogeneous spaces of ${\bf Sp}_8$ are also generated by exceptional
collections. Indeed, such a homogeneous space is a variety of partial
isotropic flags in a symplectic vector space $W$ of dimension 8, 
and is obtained by giving an ordered subset $S$ of the ordered set
$(1,2,3,4)$ (there are 16 such varieties). Let us prove that the
category $\Dd ^{b}({\rm IF}_{S})$, where the subset $S$ contains the element 1, 
has an exceptional collection. Indeed, such a homogeneous space ${\rm IF}_{S}$ is mapped onto $\Pp (W)$. 
Denote ${\rm IF}_{W}$ the flag variety and consider first the
homogeneous spaces ${\rm IF}_{1,3}$ and ${\rm IF}_{1,4}$. Let 
${\Uu}_{i}$, where $i = 1,\dots,4$, be the sequence of tautological bundles on
${\rm IF}_{W}$ as in Section 3:
\begin{equation}
0\subset {\Uu}_1\subset {\Uu}_2\subset \dots \subset {\Uu}_4 =
{\Uu}_4^{\perp}\subset {\Uu}_{3}^{\perp}\subset \dots \subset {\Uu}_1^{\perp}\subset
V\otimes \Oo _{{\rm IF}_{W}}.
\end{equation}

Denote $\pi _{1,3}$ and $\pi _{1,4}$ the projections from ${\rm IF}_{1,3}$ and ${\rm IF}_{1,4}$
onto ${\rm IF}_{1} = \Pp (W)$, respectively. The projections $\pi
_{1,3}$ and $\pi _{1,4}$ are smooth morphisms with fibers 
isomorphic to ${\rm IGr}_{2,6}$ and ${\rm LG}_{3,6}$, respectively; in
other words, these fibers being homogeneous spaces of ${\bf Sp}_6$. There is the
universal symplectic bundle ${\Uu}_{1}^{\perp}/{\Uu}_1$ on both spaces
${\rm IF}_{1,3}$ and ${\rm IF}_{1,4}$. Applying Theorem
\ref{th:Samokh1}, \cite{Sam1}, and Theorem
\ref{th:isotrgrass2-6} one obtains complete exceptional collections in
the categories $\Dd ^{b}({\rm IF}_{1,3})$ and $\Dd ^{b}({\rm IF}_{1,4})$. Any
other ${\rm IF}_{S}$ as above is either a $\Pp ^{1}$-bundle, or
a two-step iteration of $\Pp ^{1}$-bundles over either
${\rm IF}_{1,3}$, or ${\rm IF}_{1,4}$. One therefore obtains
complete exceptional collections in these cases as well. Details are left to the reader.

\begin{remark}
{\rm For symplectic groups of arbitrary ranks the above
arguments give that the study of derived categories of homogeneous
spaces of these groups can be reduced, in a way, to that of homogeneous spaces
corresponding to maximal parabolic subgroups. Indeed, if ${\bf
  P}\subset {\bf P}_{m}$ are two parabolic
subgroups of ${\bf Sp}_{2n}$, and ${\bf P}_{m}$ is maximal parabolic, then the parabolic
induction permits to map a homogeneous space ${\bf Sp}_{2n}/{\bf P}$
onto ${\bf Sp}_{2n}/{\bf P}_{m}$ with fibers being homogeneous spaces
of either symplectic or linear groups of smaller rank. Thus, if one knows the derived categories of
homogeneous spaces of groups of smaller rank and those of homogeneous
spaces that correspond to maximal parabolic subgroups of ${\bf Sp}_{2n}$ then Theorem
\ref{th:Samokh1} furnishes semiorthogonal decompositions in the
derived categories $\Dd ^{b}({\bf Sp}_{2n}/{\bf P})$. The same argument
is valid for homogeneous spaces of orthogonal groups. Finally, let us
remark that the construction of exceptional collections outlined in the proof of
Theorem \ref{th:isotrgrass2-6} may be extended to a larger
class of homogeneous spaces. Details of the computations will be given elsewhere. }
\end{remark}

\subsection*{Acknowledgements}
This paper owes its existence to numerous illuminating
conversations with A. Kuznetsov. He also made a number of valuable
remarks on a preliminary version of this text. I thank him heartily. I
would like to express also my gratitude to the University Paris 13,
where a part of this work was done.

\begin{small}

\end{small}


\begin{thebibliography}{99}
\bibitem{Be}
A. Beilinson, ``Coherent sheaves on ${\Bbb P}^n$ and Problems of
Linear Algebra'', \ Funct. Anal. Appl., 12 (1978), pp. 68--69

\bibitem{BMR}
R. Bezrukavnikov, I. Mirkovi$\acute c$  and D. Rumynin, ``Localization of
modules for a semisimple Lie algebra in prime characteristic'', \ arXiv:math.RT/0205144

\bibitem{Bo}
A. Bondal, ``Representations of associative algebras and coherent
sheaves'', (Russian) Izv. Akad. Nauk SSSR Ser. Mat. 53 (1989), no. 1,
pp. 25--44; translation in Math. USSR Izv. 34 (1990), no. 1,
pp. 23--42

\bibitem{BLL}
A. Bondal, M. Lartsen and V. Lunts, ``Grothendieck ring of
pretriangulated categories'', Int. Math. Res. Not., (2004), no. 29, pp. 1461--1495

\bibitem{BO}
A. Bondal and D. Orlov, ``Derived categories of coherent
sheaves'', Proceedings of the International Congress of Mathematicians,
Vol. II (Beijing, 2002) (Beijing), Higher Ed. Press, 2002, pp. 47--56

\bibitem{Br1}
T. Bridgeland, ``Equivalences of triangulated categories and
Fourier--Mukai transforms'', Bull. London Math. Soc. v. 31 (1999),
no. 1, pp. 25--34

\bibitem{Fu}
W. Fulton, ``Intersection Theory'',  Springer-Verlag, Berlin, New York, 1984 

\bibitem{Gro}
A. Grothendieck, ``Le groupe de Brauer'', Dix expos\'es sur la
cohomologie des sch\'emas, (1968), pp. 46--198, North-Holland, Amsterdam

\bibitem{Kap}
M. Kapranov, ``On the derived category of coherent sheaves on some
homogeneous spaces'',  Inv. Math., 92 (1988), pp. 479--508

\bibitem{Ko}
M. Kontsevich, ``Homological algebra of mirror symmetry'', Proceedings
of the International Congress of Mathematicians, Vol. 1,2 (Z\"urich,
1994) (Basel), Birkh\"auser, 1995, pp. 120--139

\bibitem{Kuznetsov1}
A. Kuznetsov, ``Hyperplane sections and derived categories'',
(Russian) Izv. RAN Ser.Mat. 70 (2006), no. 3, pp. 23--128; translation
in Izv. Math. 70 (2006), no. 3, pp. 447--547

\bibitem{Kuznetsov2}
A. Kuznetsov, ``Homological projective duality'',
arXiv:math.AG/0507292 

\bibitem{Kuznetsov3}
A. Kuznetsov, ``Derived categories of quadric fibrations and
intersections of quadrics'', arXiv:math.AG/0510670

\bibitem{Kuznetsov4}
A. Kuznetsov, ``Exceptional collections for Grassmannians
of isotropic lines'', arXiv:math.AG/0512013

\bibitem{LMS}
V. Lakshmibai, C. Musili and C.S.Seshadri, ``Cohomology of line bundles
on $\rm G/B$'', \ Annales Scientifiques de l'\'E.N.S., $4^{e}$
s\'erie, tome 7, no. 1 (1974), pp. 89 -- 137

\bibitem{Muk}
S. Mukai, ``Biregular classification of Fano threefolds and Fano
manifolds of coindex 3'', Proc. Nat. Acad. Sci., USA, 86 (1989),
pp. 3000--3002 

\bibitem{Orlov}
D. Orlov, ``Projective bundles, monoidal transformations, and derived
categories of coherent sheaves'', (Russian) Izv. Akad. Nauk SSSR Ser.Mat. 56
(1992), pp. 852--862; translation in Math. USSR Izv. 38 (1993), pp. 133--141

\bibitem{OSS}
C. Okonek, M. Schneider and H. Spindler, ``Vector bundles on
complex projective spaces'', Birkh\"auser, Boston, 1980

\bibitem{Sam1}
A. Samokhin, ``The derived category of coherent sheaves on $LG_3^{\bf C}$'',\
(Russian) Uspekhi Math.Nauk, 56, 3, (2001), pp. 592--594;
translation in Russ. Math. Surv. 56 (2001), 3, pp. 592--594

\bibitem{Sam3}
A. Samokhin, ``Tilting bundles on some Fano varieties via the Frobenius
morphism'', in preparation

\bibitem{Sch}
S. Schroer, ``Topological methods for complex-analytic Brauer
groups'', Topology, 44, (2005), pp. 875-894

\end{thebibliography}
\end{document}